\documentclass[a4paper,12pt]{article}
\usepackage{amssymb}
\input xy
\xyoption{all}
\def\s{\sigma}
\def\id{\mbox{id}}
\def\tp{\otimes}
\def\rar{\rightarrow}
\def\imp{\Rightarrow}
\def\d{\Delta}
\def\e{\varepsilon}
\def\kbar{\overline{K}}
\def\bK{\overline{K}}
\def\beq{\begin{equation}}
\def\eeq{\end{equation}}
\def\bea{\begin{eqnarray}}
\def\eea{\end{eqnarray}}
\def\nn{\nonumber}
\def\Uq2{U_q[sl_2]}
\def\Uqn{U_q[sl_n]}
\def\sl{\mathfrak{sl}}
\def\osp{\mathfrak{osp}}
\def\qqi{q-q^{-1}}

\def\CMP{\textit{Comm. Math. Phys.} }
\def\JA{\textit{J. Algebra} }
\def\lhs{\mbox{l.h.s.}}
\def\rhs{\mbox{r.h.s.}}

\begin{document}
\setlength{\baselineskip}{18pt}

\begin{titlepage}
\vspace*{2cm}
\begin{center}
{\Large\bf Weak Hopf algebras corresponding to $U_q[sl_n]$}

\vspace{1cm}
{\large
N. Aizawa${}^{\dag}$ and P. S. Isaac${}^{\ddag}$

\vspace{7mm}
${}^{\dag}$ Department of Applied Mathematics,

Osaka Women's University,
Sakai, Osaka 590-0035, JAPAN

\bigskip
${}^{\ddag}$ Graduate School of Mathematical Sciences,
University of Toyko,
3-8-1 Komaba, Meguro, Tokyo 153-8914, JAPAN
}
\end{center}
\vfill
\begin{abstract}
We investigate the weak Hopf algebras of Li based
on $U_q[sl_n]$ and Sweedler's finite dimensional example. We
give weak Hopf algebra isomorphisms between the weak generalisations of
$U_q[sl_n]$ which are ``upgraded'' automorphisms of $U_q[sl_n]$ 
and hence give a classification of these structures as
weak Hopf algebras. We also show how to decompose these examples into a
direct sum which leads to unexpected isomorphisms between their
algebraic structure.
\end{abstract}
\end{titlepage}
%
%
%
%
\setcounter{equation}{0}
\section{Introduction}

Since the introduction of quantum groups \cite{Dri1}, the importance of 
Hopf algebras has been widely recognised in both mathematics and physics. 
Generalisations of Hopf algebras have been considered,
usually motivated by some application in mathematical physics. 
The most well known example of the generalisations may be quasi-Hopf algebras 
where the coassociativity of a Hopf algebra is relaxed \cite{Dri2},
but at the same time keeping the category of modules monoidal. A similar type
of relaxation of coassociativity is also found in 
truncated quasi-Hopf algebras \cite{MS} and rational Hopf algebras \cite{Vec}. This
is where the notion of a weak co-product was introduced, such that 
$\Delta(1) \neq 1 \otimes 1,$ and was motivated by the study of
symmetries in low dimensional quantum field theory. One problem
that arose was the fact that the dual of these structures was
not associative, which led to further problems in defining
crossed products and a double construction \cite{BS}. Although these issues have 
already been addressed in
\cite{HN1} and \cite{HN2}, the question still arose as to the
possibility of defining a structure which could still provide
non-integral dimensions for the quantum field theories in a
similar way to the
weak quasi-Hopf algebras \cite{MS}, but at the
same time be coassociative.
This was the motivation behind defining the weak Hopf algebras of \cite{BS,BNS,N}.
Since these are not bialgebras, but {\it almost} bialgebras \cite{Li1},
there were also axioms required to define a weak antipode, differing
slightly from the usual ones of a Hopf algebra such
that the category of finite dimensional modules was still 
monoidal, and also with a rigidity structure defined through a
generalised antipode \cite{N}. Another question then arose as to the
possibility of defining a weak antipode on bialgebras.
Li has introduced the notion of a weak Hopf algebra to mean a
bialgebra on which is defined such a weak
antipode \cite{Li1,Li2}. In this paper, we investigate these weak
Hopf algebras as defined by Li. 

The concept of a weak Hopf algebra is rather new. 
There are not many examples known, though 
several have appeared in the literature. One example is
given by the semigroup algebra of any regular monoid which gives
a generalisation of the well known group algebra \cite{Li2}. The
other known examples were given in \cite{LD} where the authors presented
two weak Hopf generalisations of the quantised enveloping algebra $\Uq2$. 

The purposes of this paper is two fold. 
First, we wish to propose some minor adjustments to the 
examples given in \cite{LD}. Second, we extend the construction to the case of 
other known Hopf algebras such as $\Uqn$ \cite{J} and Sweedler's Hopf algebra \cite{Sw,CP}.
It is also evident that we can define weak extensions of 
quantum superalgebras in a similar way. As a consequence, we shall have a lot of new nontrivial 
examples of weak Hopf algebras. We believe that for a deeper understanding
of weak Hopf algebras as well having some insight into their applications, 
it is important to have various examples. 

The paper is organised as follows. In section \ref{wha} we give
a brief summary of the definition of weak Hopf algebra. Following
that, section \ref{wuq2} has a closer look at the examples given
in \cite{LD} and we propose slightly modified versions of these
examples. We then realise in section \ref{wuqn} that for $\Uqn$ there is a
plethora of examples which leads us to finding isomorphisms
between structures, thus giving a classification in some sense of
weak Hopf algebras corresponding to $\Uqn$. The following section
looks at weak extensions of Sweedler's famous finite dimensional
Hopf algebra, where we also show that with our construction, in
general we can decompose the weak Hopf algebra into a direct sum
of the original bialgebra with some other subalgebra. This leads
to ``unexpected'' algebra isomorphisms between structures in some cases.

%
%
%
%
\setcounter{equation}{0}
\section{Weak Hopf algebras} \label{wha}

For the reader's convenience we recall the definition of a
weak Hopf algebra in the sense of Li and Duplij \cite{Li1,LD}. 
Let $(H,\d,\e,m,u)$ be a bialgebra over a field $K$, where $\d:H\rar
H\tp H$ is the co-product, $\e:H\rar K$ is the co-unit, $m:H\tp H\rar H$
the product and $u:K\rar H$ the unit of $H$. 
The following properties define $H$:
\bea
m(m\tp\id)                        &  =  & m(\id\tp m),   \nn\\
m(u\tp\id)=                       & \id & =m(\id\tp u),  \nn\\ 
(\id\tp\d)\d                      &  =  & (\d\tp\id)\d,  \nn\\
(\e\tp\id)\d =                    & \id & =(\id\tp\e)\d, \nn\\
(m\tp m)(\id\tp\s\tp\id)(\d\tp\d) &  =  & \d\circ m,     \nn\\
\e\tp\e                           &  =  & \e\circ m.     \nn
\eea
Here $\s:H\tp H\rar H\tp H$ is the flip operator $\s(h_1\tp
h_2)=h_2\tp h_1$ for all $h_1, h_2\in H$. 
 
$H$ is a weak Hopf algebra if there is a weak antipode $T:H\rar H$
which is an algebra homomorphism satisfying the two conditions
\bea
T*\id*T   & = & T, \label{wa1} \\
\id*T*\id & = & \id \label{wa2}, 
\eea
with the convolution product $*$ defined over maps on $H$ by
$$
a*b \equiv m(a\tp b)\d :H\rar H.
$$

Note that the antipode of a Hopf algebra is a weak antipode due
to the fact that $u\circ\e :H\rar H$ is the identity of the
convolution product $*$. Recall that $S:H\rar H$ is an antipode
if it satisfies
\bea
S*\id & = & u\circ\e, \label{antipode1} \\
\id*S & = & u\circ\e. \label{antipode2}
\eea
For example, we see that
\bea
S*\id              & = & u\circ\e       \nn\\
\imp \id*S*\id     & = & \id * u\circ\e \nn\\
\imp \id * S * \id & = & \id.           \nn
\eea
In fact $S$ needs only be a left or right antipode, meaning it
satisfies only one of the two equalities (\ref{antipode1}) or
(\ref{antipode2}), in order for it to be a weak antipode.

%
%
%
\setcounter{equation}{0}
\section{Weak $\Uq2$} \label{wuq2}

In this section we give a summary of the examples of weak
$\Uq2$ presented in \cite{LD}. It is notable that the defining
relations of the ``J-weak'' quantum algebra $v\sl_q(2)$ given in that
paper can be simplified, and we give a minor adjustment to this
example and show that it is in fact a weak Hopf
algebra. We also give one other example generalising
$\Uq2$ which uses a mixture of the two examples from \cite{LD}.

We remind the reader that the usual $\Uq2$ relations to which we
refer, in terms of the four generators $E,F,K,K^{-1}$, are as follows.
\bea
K^{-1}K  & = & KK^{-1}=1,\label{r1} \\
KEK^{-1} & = & q^{2}E,\label{r2}    \\
KFK^{-1} & = & q^{-2}F,\label{r3}   \\
EF-FE    & = & \frac{K-K^{-1}}{q-q^{-1}}. \label{r4}
\eea
The co-algebra structure (co-product $\d$, co-unit $\e$) is given by 
\bea
\d(K^{\pm 1}) & = & K^{\pm 1} \tp K^{\pm 1}, \nn\\
\d(E)         & = & E\tp K + 1\tp E,         \nn\\
\d(F)         & = & F\tp 1 + K^{-1}\tp F,    \nn\\
\e(E)         & = & \e(F) = 0,               \nn\\
\e(K^{\pm 1}) & = & 1.                       \nn
\eea

It is clear that when we wish to determine the explicit action of
the antipode, we apply the definition (\ref{antipode1}) and (\ref{antipode2}) to an
arbitrary element in the algebra and solve. In all cases we can solve explicitly due to
the existence of the invertible group-like elements $1, K,
K^{-1}$. The obvious first step in generalisation to the weak
Hopf case would be to attempt to remove the invertibility of
these elements. This was the main idea in \cite{LD} when 
generalising the above definition. 

First, all weak extensions of $\Uq2$ have generators
$E,F,K,\kbar$ satisfying
\bea
K\kbar     & = & \kbar K \equiv J,\label{kk} \\
K\kbar K  & = & K,                             \ 
                \kbar K\kbar=\kbar,\label{kkk} \\
EF-FE  & = & \frac{K-\kbar}{q-q^{-1}}.  
\label{ef}
\eea
In what follows we usually write the generators with subscripts
(following \cite{LD}) to differentiate the definitions. 

{\bf Definition 1} (from \cite{LD}):  
$w\sl_{q}(2)$ is the algebra generated
by the four elements $E_w$, $F_w$, $K_w$, $\kbar_w$ satisfying
(\ref{kk},\ref{kkk},\ref{ef}) along with the relations:
\bea
K_wE_w         & = & q^{2}E_wK_w, \label{ke}          \\
\kbar_wE_w     & = & q^{-2} E_w\kbar_w,\label{kbe}    \\
K_wF_w         & = & q^{-2}F_wK_w,\label{kf}          \\
\kbar_wF_w     & = & q^{2} F_w\kbar_w,\label{kbf}
\eea
Here the invertibility of $K$ and $\kbar$ has been
relaxed, and instead of the identity, the element $J_w$ has been
introduced. It can be seen that this element $J_w$ satisfies 
$$
aJ_w = J_wa, \ \forall a\in w\sl_q(2).
$$
To demonstrate this we check 
\bea
E_w J_w & \stackrel{(\ref{kk})}{=}  & E_w K_w \kbar_w        \nn\\
        & \stackrel{(\ref{ke})}{=}  & q^{-2} K_w E_w \kbar_w \nn\\
        & \stackrel{(\ref{kbe})}{=} & K_w \kbar_w E_w    \nn\\
        & \stackrel{(\ref{kk})}{=}  & J_w E_w.
\eea
A similar calculation is performed for $F_w$ and the calculations
for $K_w$ and $\kbar_w$ are trivial. 

Also note that due to the relations (\ref{kkk}), $J_w$ is an
idempotent. Namely,
$$
J_w^2=J_w.
$$
The co-algebra structure is defined as follows. 
The co-product and co-unit are respectively given by
\bea
\d_w(E_w)     & = & 1\tp E_w + E_w\tp K_w,     \nn\\
\d_w(F_w)     & = & F_w\tp 1 + \kbar_w\tp F_w, \nn\\
\d_w(K_w)     & = & K_w\tp K_w,                \nn\\
\d_w(\kbar_w) & = & \kbar_w\tp \kbar_w,        \nn\\
\e_w(E_w)     & = & \e_w(F_w) = 0,               \nn\\
\e_w(K_w)     & = & \e_w(\kbar_w) = 1.           \nn
\eea 
It can be verified that they are both algebra homomorphisms so
that 
$$
\d_w(xy)=\d_w(x)\d_w(y)
$$ 
and 
$$
\e_w(xy)=\e_w(x)\e_w(y)
$$ for all $x,y\in w\sl_{q}(2)$, thus preserving the defining
relations. 
With this co-product a corresponding weak antipode 
can be determined by solving equations (\ref{wa1})
and (\ref{wa2}) with the above co-product. The only possible
weak antipode in this case is
\bea
T_w(1)       & = & 1,           \nn\\
T_w(K_w)     & = & \kbar_w,     \nn\\
T_w(\kbar_w) & = & K_w,         \nn\\
T_w(E_w)     & = & -E_w\kbar_w, \nn\\
T_w(F_w)     & = & -K_wF_w.     \label{wasl2}
\eea
It can be shown that $T_w$ is an algebra anti-homomorphism,
that is, $T(ab)=T(b)T(a)$. Note that with the above bialgebra 
structure it is not possible
to determine an antipode in the usual sense. As we mentioned
previously, this was the motivation for relaxing (\ref{r1}) to
(\ref{kk}) in order to provide weak antipodes which are not
antipodes. For example, to solve the equation
$$
S*\id (K) = \e(K)1 \ \imp \ S(K)K=1,
$$
we would need an inverse of the element $K$. 

Another possible definition given in \cite{LD} is the
following. 

{\bf Definition 2} (from \cite{LD}):
$v\sl_q(2)$ is the algebra generated
by the four elements $E_v$, $F_v$, $K_v$, $\kbar_v$ satisfying
(\ref{kk},\ref{kkk},\ref{ef}) along with the relations:
\bea
K_vE_v\kbar_v       & = & q^{2}E_v,\label{ke2}              \\
K_vF_v\kbar_v       & = & q^{-2}F_v.\label{kf2} 
\eea
In this case, $J_v=K_v\kbar_v$ satisfies the relation
\beq
J_v a = a J_v = a, 
\label{ja}
\eeq
for $a=E_v,\ F_v,\ K_v,\ \kbar_v$ (and hence $J_v$). 
To demonstrate,
we have
\bea
E_v J_v & \stackrel{(\ref{kk})}{=} & E_vK_v\kbar_v \nn\\
& \stackrel{(\ref{ke2})}{=} & q^{-2}K_vE_v\kbar_v K_v \kbar_v \nn\\
& \stackrel{(\ref{kkk})}{=} & q^{-2}K_vE_v\kbar_v
\ (\stackrel{(\ref{ke2})}{=} E_v) \nn\\
& \stackrel{(\ref{kkk})}{=} & q^{-2}K_v\kbar_v K_vE_v\kbar_v \nn\\
& \stackrel{(\ref{ke2})}{=} & K_v\kbar_v E_v \nn\\
& \stackrel{(\ref{kk})}{=} & J_v E_v.  \nn
\eea
For the generator $F_v$, a similar calculation can be done. For
the cases $K_v$ and $\kbar_v$, the calculation is trivial. 
The most remarkable consequence of this property is that the
analogue of
defining relation (\ref{ef}) presented in {\cite{LD} which was in the
form
$$
E_vJ_vF_v - F_vJ_vE_v = \frac{K_v-\kbar_v}{q-q^{-1}},
$$
reduces to (\ref{ef}) by the above argument. Therefore in what
follows we shall
always use relation (\ref{ef}) and not the relation above.

The co-algebra structure for this second definition is as
follows. 
\bea
\d_v(E_v) & = & J_v\tp E_v + E_v\tp K_v, \nn\\
\d_v(F_v) & = & F_v\tp J_v + \kbar_v \tp F_v, \nn\\
\eea
with the remaining actions coinciding precisely with the
case of definition 1.

Moreover, relations (\ref{ke2}) and (\ref{kf2}) can be manipulated 
to those of definition 1. We demonstrate that
\bea
K_vE_v & \stackrel{(\ref{kkk})}{=} & K_v\kbar_vK_vE_v \nn\\
       & \stackrel{(\ref{kk}),(\ref{ja})}{=} & K_vE_v\kbar_vK_v \nn\\
       & \stackrel{(\ref{ke2})}{=} & q^2E_vK_v.  \nn
\eea
The other relations can be verified in a similar way. Although
the co-product is different to that of definition 1, there
exists a weak antipode which is the same as (\ref{wasl2}) in the case of
definition 1.

This indicates that much of the discussion in \cite{LD} relating
to $v\sl_q(2)$ is redundant. However, we would like to make it
clear that we consider the paper \cite{LD} rich in ideas and an
inspiration to our current investigations. 

There are other possibilities for defining weak
extensions of $\Uq2$. These involve mixtures of definition 1 and
definition 2 over the generators $E,F$. For example, we can
say that one case is where $E$ satisfies the relations (\ref{ke}), 
(\ref{kbe}) and $F$ satisfies (\ref{kf2}),
along with all the other relations common to both definition 1
and definition 2. The co-product would then have the action
\bea
\d(E) & = & 1\tp E + E\tp K, \nn\\
\d(F) & = & F\tp J + \kbar\tp F, \nn
\eea
along with the usual group-like co-product for $K$ and $\kbar$.
The weak antipode would still be the same as in definitions 1 and
2. 

We can also swap this mixture of definitions and say that $E$
satisfies those relations of definition 2, but $F$ satisfies the
relations of definition 1. This case is actually isomorphic to
the first mixture, as we shall see later. In the section on Weak
$\Uqn$ we give a more formal way of notating such mixtures.

So we now have some clues as to how we may approach the problem
of defining weak extensions of $\Uqn$. It is clear that there
will be many possible combinations of generators satisfying
either of the two definitions in the
general case. This then begs the question: how would we know
which mixtures of the two definitions lead to isomorphic algebras?
To this end we have an important observation regarding some of the
automorphisms of the original quantum algebra $\Uqn$ which
``lift up'' to isomorphisms between weak Hopf structures. We
shall look at these isomorphisms in more detail in the next section.

In general, we say that a generator satisfying the 
relations of definition 1 is of type 1, and is type 2 if it
satisfies the relations of definition 2.
%
%
%
%
\setcounter{equation}{0}
\section{Weak $\Uqn$} \label{wuqn}

\subsection{Mixing definitions}

For the case of $w\sl_q(n)$, which 
has simple generators $E_i$, $F_i$, $K_i$ and $\kbar_i$ ($i=1,\ldots ,n-1$), we
can choose either definition 1 or 2 to describe the relations
between any $E_i$ and the $K_j/\kbar_j$ and similarly for any $F_i$. This is what is meant
by the word ``mixture''. 
The relations satisfied by the generators are as follows, for all
$i,j$ unless specified otherwise;
\bea
K_iK_j = K_jK_i,\ \kbar_i\kbar_j = \kbar_j\kbar_i, \ K_i\kbar_j = \kbar_j K_i,\ K_i\kbar_i & = & J,\nn\\
JK_j = K_jJ = K_j, \ J\kbar_j = \kbar_jJ & = & \kbar_j,\nn\\
E_iF_j - F_jE_i & = & \delta_{ij}\frac{K_i-\kbar_i}{q-q^{-1}},
\nn\\
E_i^2E_{i\pm 1}-(q+q^{-1})E_iE_{i\pm 1}E_i+E_{i\pm 1}E_i^2 & = &
0, \nn\\
F_i^2F_{i\pm 1}-(q+q^{-1})F_iF_{i\pm 1}F_i+F_{i\pm 1}F_i^2 & = &
0, \nn\\
E_iE_j  =  E_jE_i, \ F_iF_j & = & F_jF_i, \ |i-j|\geq 2 \label{slnrels}.
\eea
We also need to specify the relations between the $E_i$ and the
$K_j$ for example. Let $a_{ij}$ denote the Cartan matrix for
$sl(n)$, $a_{ii}=2$, $a_{i,i\pm 1} = -1$ and zero otherwise. If $E_i$ satisfies
\beq
K_jE_i = q^{a_{ij}}E_iK_j,\ E_i\kbar_j =q^{a_{ij}}\kbar_jE_i,\
\forall j,  \label{type1}
\eeq
we say $E_i$ satisfies definition 1, or simply $E_i$ is type 1.
However, if $E_i$ satisfies 
\beq
K_jE_i\kbar_j = q^{a_{ij}}E_i, \ \forall j, \label{type2}
\eeq
we say $E_i$ satisfies definition 2, or simply $E_i$ is type 2.
The same convention holds for $F_i$ by replacing $E_i$ with
$F_i$ and $a_{ij}$ with $-a_{ij}$ in the above relation. Notice
also that $J$ is defined for all $i$, so for example
$$
J=K_i\kbar_i=K_j\kbar_j,\ \ i\neq j.
$$

The co-product has the following action;
\bea
\d(K_i) & = & K_i\tp K_i, \nn\\
\d(\kbar_i) & = & \kbar_i\tp\kbar_i,\nn\\
\d(E_i) & = & \left\{ 
            \begin{array}{ll}
               1\tp E_i+E_i\tp K_i, & \mbox{$E_i$ is type 1} \\
               J\tp E_i+E_i\tp K_i, & \mbox{$E_i$ is type 2} 
            \end{array} 
              \right. \nn\\
\d(F_i) & = & \left\{ 
         \begin{array}{ll}
           F_i\tp 1+\kbar_i\tp F_i, & \mbox{$F_i$ is type 1} \\
           F_i\tp J+\kbar_i\tp F_i, & \mbox{$F_i$ is type 2}
         \end{array}
              \right. \label{uqncoprod}
\eea
while the action of the co-unit is
$$
\e(1)=\e(K_i)=\e(\kbar_i)=1,\ \ \e(E_i)=\e(F_i) = 0.
$$
The weak antipode $T$ will always have the form
\bea
T(1)       & = & 1,           \nn\\
T(K_i)     & = & \kbar_i,     \nn\\
T(\kbar_i) & = & K_i,         \nn\\
T(E_i)     & = & -E_i\kbar_i, \nn\\
T(F_i)     & = & -K_iF_i,     \nn
\eea
regardless of the type of the generators $E_i$ and $F_i$.

In order to notate these mixtures for $w\sl_q(n)$ we use a binary
notation, where a 1 indicates the use of a type 1 
generator and a 0 indicates the use of a type 2 generator. We list the
$2(n-1)$ simple generators $E_i$ and $F_i$, starting with the $E_i$ followed by
the $F_i$. We then write down a list of $0$'s and $1$'s in the order
corresponding to the generators determined by their 
type. 
This then gives an integer from $0$ to
$2^{2(n-1)}-1$ in binary representation which contains all the information as to
which particular mixture of definition we are using for the
relations between the generators $E_i$ and $F_i$ and all the $K_j/\kbar_j$. 
We denote this integer $d$ and the algebra is expressed as $w\sl^d_q(n)$.
In total there are $2^{2(n-1)}$ possible mixtures for
$w\sl^d_q(n)$.

Note that we cannot have different definitions for the relations
between the same generator with different $K_i$'s because 
the co-product could not possibly be consistent with those defining
relations.

For example, in the case of $w\sl_q(4)$ we have the simple
generators (not including the $K_i$), $E_1, E_2, E_3, F_1, F_2, F_3$. 
Hence there are
$2^6=64$ different possibilities for relations with the $K_i$. 
The notation $w\sl^{43}_q(4)$
has the following meaning. Since the number $43$ has the binary
representation $101011$, this is interpreted to mean that the
simple generators $E_1, E_3, F_2, F_3$ are type 1 with 
the remaining ones $E_2, F_1$ being type 2. 
This information is determined by
superimposing the list of binary digits $\{1,0,1,0,1,1\}$ with the list 
of simple generators in the
order $\{E_1, E_2, E_3, F_1, F_2, F_3\}$.

It should also be noted that the algebra $w\sl^3_q(2)$ coincides
with $w\sl_q(2)$ given in \cite{LD} (and in section \ref{wuq2}
above) and the example $v\sl_q(2)$ of \cite{LD} is precisely $w\sl^0_q(2)$ in
our notation.

\subsection{Isomorphic structures}

Now we look in more detail at the weak Hopf algebras of type
$\Uqn$ using mixtures of the two types of generators. In some cases
where $d_1\neq d_2$, there exists a weak Hopf algebra isomorphism 
$w\sl^{d_1}_q(n)\simeq w\sl^{d_2}_q(n)$. It is therefore worth
investigating all possible isomorphisms in order to classify the
weak extensions based on our criteria. As we shall see in this
section, the isomorphisms are derived from a subset of the set of
automorphisms on the algebra $\Uqn$. In other words, a subset of
the automorphisms on $\Uqn$ ``lift up'' to isomorphisms between
the weak Hopf extensions. The reason only a subset can be
considered, as we
shall see later, is because some of the automorphisms of $\Uqn$
lose their invertibility when upgraded to act on the weak $w\sl^d_q(n)$,
so therefore cannot be isomorphisms.

If $(A,\d,\e,T)$ and $(B,\d',\e', T')$ are weak Hopf algebras, then a weak Hopf algebra 
isomorphism $\psi:A\rar B$ is an invertible algebra homomorphism satisfying
\bea
(\psi\tp\psi)\circ\d & = &\d'\circ\psi, \label{consistency} \\
\e & = & \e'\circ\psi, \label{consistency2} \\
\psi\circ T & = & T'\circ\psi. \label{consistency3}
\eea

For example, consider the algebra $w\sl^1_q(2)$. This has generators
$E^{(1)},F^{(1)},K^{(1)},\kbar^{(1)}$ (and $1^{(1)}$) satisfying
\bea 
K^{(1)}\kbar^{(1)}   & = & \kbar^{(1)} K^{(1)} \equiv  J^{(1)}, \nn\\
K^{(1)}\kbar^{(1)} K^{(1)} & = & K^{(1)}, \ \ \kbar^{(1)} K^{(1)}\kbar^{(1)} = \kbar^{(1)}, \nn\\
K^{(1)}F^{(1)}       & = & q^{-2}F^{(1)}K^{(1)}, \ \ \kbar^{(1)} F^{(1)} = q^2F^{(1)}\kbar^{(1)},\nn\\
K^{(1)}E^{(1)}\kbar^{(1)}  & = & q^2E^{(1)},\nn\\
E^{(1)}F^{(1)}-F^{(1)}E^{(1)}    & = & \frac{K^{(1)}-\kbar^{(1)}}{q-q^{-1}},\nn
\eea
since the binary representation of $1=\{0,1\}$ is superimposed with
the list of generators $\{E^{(1)}, F^{(1)}\}$ and so $E^{(1)}$ is type 2
and $F^{(1)}$ is type 1.
Let us now consider the co-algebra structure of this algebra.
The co-product $\d$ and co-unit $\e$ are respectively given by
\bea
\d(K^{(1)})     & = & K^{(1)}\tp K^{(1)}, \nn\\
\d(\kbar^{(1)}) & = & \kbar^{(1)}\tp\kbar^{(1)},\nn\\
\d(E^{(1)})     & = & J^{(1)}\tp E^{(1)} + E^{(1)}\tp K^{(1)}, \nn\\
\d(F^{(1)})     & = & F^{(1)}\tp 1^{(1)} + \kbar^{(1)}\tp F^{(1)}, \nn\\
\e(E^{(1)})     & = & \e(F^{(1)}) = 0,\nn\\
\e(K^{(1)})     & = & \e(\kbar^{(1)}) = 1.\nn
\eea
Now consider the map $\psi: w\sl^1_q(2) \rar w\sl^2_q(2)$ defined
by the action
\bea
\psi(E^{(1)}) & = & F^{(2)}, \nn\\
\psi(F^{(1)}) & = & E^{(2)},\nn\\
\psi(K^{(1)}) & = & \kbar^{(2)},\nn\\
\psi(\kbar^{(1)}) & = & K^{(2)},\nn
\eea
where we have employed an obvious notation with superscripts.
This map derives from the so-called Cartan involution on $\Uq2$.
In the weak case it can be seen to be a weak Hopf algebra isomorphism since
it preserves the generator type (that is, it is consistent with the defining
relations) and is also consistent with equations
(\ref{consistency})-(\ref{consistency3}).

In general the rule is that such
an isomorphism must map a type 1 generator into a type 1
generator and similarly for type 2. 
We demonstrate the sort of calculation required to show
consistency with the relations. Take for example
\bea
\lhs = \psi(K^{(1)})\psi(E^{(1)})\psi(\kbar^{(1)}) & = & \kbar^{(2)}F^{(2)}K^{(2)}
\nn\\
& = & q^2 \kbar^{(2)}(K^{(2)}F^{(2)}\kbar^{(2)})K^{(2)}\nn\\
& = & q^2F^{(2)}\nn\\
& = & q^2\psi(E^{(1)})= \rhs \nn
\eea
The other relations can be realised in a similar fashion.
To demonstrate consistency with the co-product, we can use equation
(\ref{consistency}) to determine $\d'$ in this case. For example, applying
both sides of (\ref{consistency}) to $E^{(1)}$ gives 
\bea
(\psi\tp\psi)\d(E^{(1)}) & = & (\psi\tp\psi)(J^{(1)}\tp
E^{(1)}+E^{(1)}\tp K^{(1)}) \nn\\
& = & J^{(2)}\tp F^{(2)}+F^{(2)}\tp \kbar^{(2)}, \nn\\
\d'(\psi(E^{(1)})) & = & \d'(F^{(2)}). \nn
\eea
This then gives the action of $\d'$ on $F^{(2)}$. The remaining
actions are
\bea
\d'(E^{(2)}) & = & E^{(2)}\tp 1^{(2)}+K^{(2)}\tp E^{(2)}, \nn\\
\d'(K^{(2)}) & = & K^{(2)}\tp K^{(2)},\nn\\
\d'(\kbar^{(2)}) & = & \kbar^{(2)}\tp \kbar^{(2)}.\nn
\eea
Note that this is not the co-product given in equations
(\ref{uqncoprod}), but it is in fact the opposite co-product,
$\d^\s = \s\circ\d$ ($\s$ being the flip operator), which we know 
from the theory of bialgebras is a perfectly
acceptable one. It is also straightforward to verify
(\ref{consistency2}) holds. Because the action of the weak antipode is dependent
on the co-product, $T'$ will be different to the one presented
earlier. It is straightforward to verify that it does indeed
exist, and that equation (\ref{consistency3}) is satisfied.

Although there are undoubtedly many other possibilities to extend
$\Uq2$ to a weak structure, the extensions presented in this
paper based on that of \cite{LD} total three, namely,
$w\sl^0_q(2)$, $w\sl^1_q(2)\simeq w\sl^2_q(2)$ and $w\sl^3_q(2)$.
In each case the weak antipode has the same action, that of
(\ref{wasl2}).

\subsection{Number of unique structures}

We now address the question of the number of possible weak Hopf algebra isomorphisms $\psi:
w\sl^{d}_q(n)\rar w\sl^{d'}_q(n)$, $d\neq d'$. 

It is well known that for $\Uqn$ there are several types of
automorphisms \cite{T}. The most relevant to this paper are
the Dynkin diagram automorphisms and the Cartan
involution, since they give rise to Hopf algebra automorphisms
and anti-automorphisms respectively. 
We consider maps $\rho_d$ and $\omega_d$ which have the same
actions as the Dynkin diagram automorphism and the Cartan
involution respectively, but applied to the weak Hopf
algebra $w\sl^d_q(n)$. These maps then become isomorphisms
between weak Hopf structures. Their actions are given by
$$
\begin{array}{lcllcllcllcl}
\rho_d(E^{(d)}_i)     & = &  E^{(d')}_{n-i}, &
\rho_d(F^{(d)}_i)     & = &  F^{(d')}_{n-i}, &
\rho_d(K^{(d)}_i)     & = &  K^{(d')}_{n-i}, &
\rho_d(\kbar^{(d)}_i) & = & \kbar^{(d')}_{n-i}, \\
\omega_d(E^{(d)}_i)     & = &  F^{(d'')}_i, &
\omega_d(F^{(d)}_i)     & = &  E^{(d'')}_i, &
\omega_d(K^{(d)}_i)     & = & \kbar^{(d'')}_i, &
\omega_d(\kbar^{(d)}_i) & = & K^{(d'')}_i,
\end{array}
$$
where the indices $d$ (corresponding to the source), $d'$ and $d''$ (corresponding
to the targets) are used to
differentiate between the structures. Note that $\rho_d$ and
$\omega_d$ map into different spaces which justifies the use of the
different indices $d'$ and $d''$.
With these actions, it can be easily verified that 
$$
\rho_{d'}\circ\rho_d = \id,\ \ \omega_{d''}\circ\omega_d = \id.
$$

Lusztig \cite{Lus1} has also given a set of algebra automorphisms 
defined on the quantised enveloping algebras. However, when applied
to our weak generalisations, they are found to be 
non-invertible and therefore are not isomorphisms between weak
Hopf algebras. 

Although we know of the existence of other algebra isomorphisms which
exist between structures, the only known weak Hopf algebra isomorphisms
are $\rho_d$ and $\omega_d$. We will comment more on these
algebra isomorphisms in section \ref{swee}.

One important point is that $\rho_d$ and $\omega_d$ both preserve
the generator type, so for example if $E^{(d)}_i$ is a type 1
generator, so are $E^{(d')}_{n-i}$ and $F^{(d'')}_i$. Therefore
$\rho_d$ and $\omega_d$ must correspond to maps (say
$r_d$ and $w_d$ respectively) defined on the non-negative
integers such that
\bea
\rho_d:w\sl^d_q(n) \rar w\sl^{r_d(d)}_q(n),&&\nn\\
\omega_d:w\sl^d_q(n) \rar w\sl^{w_d(d)}_q(n).&&\nn
\eea
Once we know the action of the maps $r_d$ and $w_d$, this should
allow us to be able to determine which structures are isomorphic
and hence lead to a classification.

To this end, we write $d$ in terms of its binary expansion 
$$d=(d_0,d_1,\ldots,d_{n-2}|d_{n-1},\ldots,d_{2n-3}),$$ where the bar
separates the values representing the $E_i$ and $F_i$, and where
the $d_i$ have values of either 0 or 1. Then $r_d(d)$ and
$w_d(d)$ have the expansions
\bea
w_d(d) & = & (d_{n-1},\ldots,d_{2n-3}|d_0,\ldots,d_{n-2}), \nn\\
r_d(d) & = & (d_{n-2},\ldots,d_0|d_{2n-3},\ldots,d_{n-1}).\nn
\eea
In terms of the components of the binary expansion we have
\bea
w_d(d_k) & = & d'_{3n-4-k\bmod 2(n-1)},\nn\\
r_d(d_k) & = & d''_{n-1+k\bmod 2(n-1)}.\nn
\eea
This simplifies the problem of determining isomorphic structures and
allows us to explicitly count the number of unique structures for
each $n$.

It is worth noting that $\rho\circ\omega = \omega\circ\rho$, 
so the only isomorphisms we need to
consider are $\rho$, $\omega$ and $\rho\circ\omega$ (we have dropped the
subscripts for convenience). 
According to this prescription, there can only be at most four
structures which are isomorphic. In some cases there could be two
and in others there could be no isomorphisms. These cases are
referred to below as degenerate.
The situation can be summarised in the following diagram
$$
\xymatrix{
& r(d) \ar[r] & w\circ r(d) \\
d \ar[ur]\ar[dr] & & \\
& w(d) \ar[r] & r\circ w(d) \ar@{=}[uu]
}
$$
where in some cases the arrows could be equalities, in which case
there would be one of the afore mentioned degeneracies. In order
to count the number of unique (non-isomorphic) weak extensions of
$w\sl^d_q(n)$, we list all the possible degenerate cases;
\bea
(1) && w(d)=d,\nn\\
(2) && r(d)=d,\nn\\
(3) && r(d)=w(d),\nn
\eea
and consider their intersection $(1)\cap (2)\cap (3)$, union $(1)\cup(2)\cup(3)$ and their union's complement
$\overline{(1)\cup(2)\cup(3)}$
in order to count the total number of unique cases. We note that
some of these cases will lead to exactly two isomorphic
structures, and combinations of the above cases will lead to no
isomorphic structures. We aim to separate each of these situations and
then add the number of structures relating to each. 

\noindent Case (1): The only $d$ satisfying $w(d)=d$ is of the
form
$$
d=(d_0,d_1,\ldots,d_{n-2}|d_0,d_1,\ldots,d_{n-2}).
$$
Therefore the total number of cases satisfying case (1) is
$2^{n-1}$.

\noindent Case (2): We separate this case into two cases
corresponding to $n$ being odd and even. For $n=2m+1$ the only
$d$ satisfying case (2) is of the form
$$
d=(d_0,d_1,\ldots,d_{m-1},d_{m-1},\ldots,d_0|d_{2m},d_{2m+1},\ldots,d_{3m-1},d_{3m-1},\ldots,d_{2m}),
$$ 
so there are $2^{2m}=2^{n-1}$ possibilities. For $n=2m$ the only
$d$ satisfying case (2) has the form
$$
d=(d_0,\ldots,d_{m-2},d_{m-1},d_{m-2},\ldots,d_0|d_{2m-1},\ldots,d_{3m-3},d_{3m-2},d_{3m-3},\ldots,d_{2m-1}),
$$
so there are $2^{2(m-1)+2}=2^{2m}=2^n$ possibilities.

\noindent Case (3): The only $d$ satisfying this case is of the
form
$$
d=(d_0,d_1,\ldots,d_{n-2}|d_{n-2},\ldots,d_0),
$$
so there are $2^{n-1}$ possibilities.

\noindent (1) $\cap$ (2) $\cap$ (3): Once again we treat the case
for $n$ is even and odd separately. 
For $n=2m+1$, $d$ is of the form
$$
d=(d_0,d_1,\ldots,d_{m-1},d_{m-1},\ldots,d_0|d_0,d_1,\ldots,d_{m-1},d_{m-1},\ldots,d_0),
$$ 
so there are $2^m=2^{(n-1)/2}$ possibilities. For $n=2m$, $d$ is
of the form
$$
d=(d_0,\ldots,d_{m-2},d_{m-1},d_{m-2},\ldots,d_0|d_0,\ldots,d_{m-2},d_{m-1},d_{m-2},\ldots,d_0),
$$
giving $2^m=2^{n/2}$ possibilities.

\noindent $\overline{(1)\cup(2)\cup(3)}$: Combining the above three
cases and subtracting twice their intersection gives the union,
for which there are $3.2^{2m}-2^{m+1}$ possibilities for $n=2m+1$
and $2^{2m+1}-2^{m+1}$ possibilities for $n=2m$. The compliment
of the union then has $2^{4m}-3.2^{2m}+2^{m+1}$ possibilities for
$n=2m+1$ and $2^{4m-2}-2^{2m+1}+2^{m+1}$ possibilities for $n=2m$.

To calculate the exact number of unique structures we consider
the fact that case (1) without cases (2) or (3) (and
permutations) will have precisly 2 structures which are
isomorphic, or put another way, isomorphic structures in these
cases come in
pairs. Therefore we need to half the number obtained above when counting the
total number of structures. Similarly for the structures which
do not fall into these degenerate cases. There will be exactly 4
isomorphic structures so we need to divide the number
corresponding to $\overline{(1)\cup(2)\cup(3)}$ by 4. 

Therefore the number of non-isomorphic structures, say $Z_n$, is
\bea
Z_{2m+1} & = & \frac{2^{2m}-2^m}{2}+\frac{2^{2m}-2^m}{2}+\frac{2^{2m}-2^m}{2}+2^m
+ \frac{2^{4m}-3.2^{2m}+2^{m+1}}{4} \nn\\
         & = & 2^{4m-2}+\frac{3}{4}.2^{2m},\nn\\
Z_{2m} & = & \frac{2^{2m-1}-2^m}{2}+\frac{2^{2m}-2^m}{2}+\frac{2^{2m-1}-2^m}{2}+2^m 
+ \frac{2^{4m-2}-2^{2m+1}+2^{m+1}}{4} \nn\\
       & = & 2^{4m-4}+2^{2m-1}.\nn
\eea
Putting these two cases together gives
$$
Z_n=2^{n-4}(7+(-1)^n+2^n),
$$
which is the number of unique weak Hopf structures corresponding
to $w\sl^d_q(n).$

This formula for $Z_n$ has been verified up to $n=10$ by directly applying the
maps $\rho$, $\omega$ and $\rho\circ\omega$ and then counting the number of unique
structures. To give the reader an idea of the number of
structures, we have the table below.

$$
\begin{tabular}{|l|l|l|l|l|l|l|l|l|l|l}
\hline
$n$  & 2 & 3 & 4 & 5 & 6 & 7 & 8 & 9 & 10 \\
\hline
$Z_n$ & 3 & 7 & 24 & 76 & 288 & 1072 & 4224 & 16576 & 66048 \\
\hline
\end{tabular}
$$

We also list for $n\leq 4$ all the values of $d$, putting
isomorphic values in brackets $\{,\}$.
For $n=2$ we have already determined that the values
$$
d=0,\{1,2\}, 3
$$
give the 3 unique structures. For $n=3$ the values of $d$ for
the 7 structures are 
$$
d=0,\{1,2,4,8\},\{3,12\},\{5,10\},\{6,9\},\{7,11,13,14\},15.
$$
For $n=4$ the 24 values of $d$ are
\bea
d & = & 0,
\{1,4,8,32\},
\{2,16\},
\{3,6,24,48\},
\{5,40\},
\{7,56\}, 
\{9,36\}, 
\{10,17,20,34\}, \nn\\
& & \{11,25,38,52\}, 
\{12,33\}, 
\{13,37,41,44\}, 
\{14,28,35,49\}, 
\{15,39,57,60\}, 
18, \nn\\
& & 
\{19,22,26,50\}, 
\{21,42\}, 
\{23,58\},
\{27,54\}, 
\{29,43,46,53\}, 
\{30,51\}, \nn\\
& & 
\{31,55,59,62\}, 
45, 
\{47,61\}, 
63. \nn
\eea
All cases up to $n=10$ have been calculated, but are obviously too
unwieldy to include in the article.

%
%
%
%
\setcounter{equation}{0}
\section{Direct sum decomposition and Sweedler's example}
\label{swee}

We now look in more detail at the algebraic structure and show
that the upgraded quantised enveloping algebra automorphisms are
not the only algebra isomorphisms between the various
$w\sl^d_q(n)$. 

First recall Sweedler's example \cite{Sw} of a finite dimensional Hopf
algebra, denoted $H$. $H$ is generated by elements $I, G, X$ (where $I$
is the identity element) satisfying the relations
\bea
G^2 & = & I, \nn\\
GX & = & -XG, \nn\\
X^2 & = & 0. \nn
\eea
The co-product is given by
\bea
\d(G) & = & G\tp G, \nn\\
\d(X) & = & X\tp G+I\tp X, \nn
\eea
and the co-unit given by
$$
\e(G)=1=\e(I),\ \ \ \e(X)=0.
$$
The antipode $S$ is given by the action
$$
S(G)=G, \ \ S(I) = I, \ \ S(X)=GX.
$$
Clearly $H$ is 4 dimensional with basis $\{I, G, GX, X\}$. 

In order to give an example of a weak Hopf algebra based on this
structure with generators $\{1,g,x\}$ (we now use lower case
symbols), instead of using the relation $g^2=1$, we impose
the relation $g^3=g$. Moreover, we can choose either the relation $gx=-xg$, in
which case we refer to $x$ as a type 1 generator (analagous to
the notion discussed at the
end of section \ref{wuq2}), or we can choose the relation 
$gxg = -x$, in which case we call $x$ a type 2 generator.
A type 2 generator is also a type 1
generator, but not conversely, since $g^2\neq 1$.

For the first case, we choose $x$ to be type 1. Denote the
algebra by $H_1$. The following relations are satisfied.
\bea
g^3 & = & g, \nn\\
gx & = & -xg, \nn\\
x^2 & = & 0, \nn
\eea
along with the same co-product and co-unit as in the usual Hopf case (given above).
Solving equations (\ref{wa1}) and (\ref{wa2}) gives the weak
antipode
\bea
T(1) & = & 1, \nn\\
T(g) & = & g, \nn\\
T(x) & = & gx, \nn
\eea
which has the same action as the antipode from the Hopf case but
nevertheless is not an antipode.
The defining relations imply that $H_1$ is $6$ dimensional with
basis
$$
\{ 1,g,g^2,x,gx,g^2x \}.
$$
Note that the element $g^2$ is a central idempotent. This is
easily verified with the defining relations. 

It is with this example that we demonstrate explicitly how to
obtain a direct sum decomposition for an algebra with a central
idempotent. This procedure will then be extended to the case of
$w\sl^d_q(n)$.

$H_1$ has a direct sum decomposition
$$
H_1 = H_1^0 \oplus H_1^1,
$$
where $H_1^0$ is the subalgebra with basis $\{(1-g^2)x,1-g^2 \}$,
on which multiplication by $g^2$ is zero (indicated in the
superscript), and $H_1^1$ is the
subalgebra with basis $\{ g,g^2,gx,g^2x \}$ on which
multiplication by $g^2$ is the 
identity (also indicated in the superscript). In fact, these two
subalgebras are determined by setting
\bea
H_1^0 & = & (1-g^2)H_1, \nn\\
H_1^1 & = & g^2H_1.
\eea
It is straightforward to verify that the map $\psi:H_1^1 \rar H$
with the action
\bea
\psi(g) & = & G, \nn\\
\psi(g^2) & = & I, \nn\\
\psi(gx) & = & GX,\nn\\
\psi(g^2x) & = & X,\nn
\eea
defines a weak Hopf algebra isomorphism, 
where $I,G,X$ are the generators
of the original Sweedler Hopf algebra $H$.
Since $H$ appears as a subalgebra of $H_1$, we can simply apply
$\psi^{-1}\tp\psi^{-1}$ to the R-matrix of $H$ (see \cite{CP}) to obtain an
R-matrix ${\cal R}$ of $H_1$ satisfying
\bea
{\cal R}\d(a) & = & \s\circ\d(a){\cal R}, \ \forall a\in H_1 \nn\\
{\cal R}_{13} {\cal R}_{23} & = & (\Delta \otimes \id)({\cal R}),
\nn\\ 
{\cal R}_{13} {\cal R}_{12} & = & (\id \otimes \Delta)({\cal R}).\nn
\eea
Such an R-matrix is then given by
$$
{\cal R}=g^2\tp g^2-2p\tp p+\alpha(g^2x\tp g^2x-2g^2x\tp px+2px\tp px),
$$
where $p=(g^2-g)/2$ and $\alpha$ is an arbitrary parameter.
This ${\cal R}$ is not invertible, but it satisfies the
regularity condition \cite{LD}
\bea
{\cal R}\hat{\cal R} {\cal R} & = & {\cal R}, \\
\hat{\cal R}{\cal R}\hat{\cal R} & = & \hat{\cal R},
\eea
where 
$$
\hat{\cal R}=g^2\tp g^2-2p\tp p+\alpha(g^2x\tp g^2x-2px\tp g^2x+2px\tp px).
$$

It should also be noted that the Sweedler Hopf algebra $H$ also
appears as a subalgebra of $H_1$ with
its generators defined by
\bea
I & = & 1, \nn\\
G & = & 1+g-g^2, \nn\\
X & = & (1+\alpha g)gx,\nn
\eea
where $\alpha$ is an an arbitrary constant. However, this is just 
an observation and has no consequence to the results of our paper,
since this subalgebra is only isomorphic to $H$ as an algebra,
not a bialgebra.

Now we look at the algebra $H_2$, which corresponds to the choice
of the generator $x$ to be of type 2. This implies the following
relations;
\bea
g^3 & = & g, \nn\\
gxg & = & -x, \nn\\
x^2 & = & 0. \nn
\eea
The only difference with the co-product in this case is with the
action defined on the generator $x$, which is now given by
$$
\d(x) = x\tp g + g^2\tp x,
$$
and the co-unit is the same as usual.
The algebra $H_2$ is $5$ dimensional with basis $\{ 1, g, x, gx,
g^2 \}$. Note that $g^2$ is a central idempotent. Defining 
\bea
H_2^0 & = & (1-g^2)H_2, \nn\\
H_2^1 & = & g^2H_2,\nn
\eea
the decomposition
$$
H_2 = H_2^0 \oplus H_2^1
$$
still holds, where the superscripts still refer to the action of $g^2$, but
now $H_2^0$ has basis $\{ 1-g^2 \}$ and $H_2^1$ has basis $\{
g^2, g, x, gx \}$.

Once again it is possible to verify that there exists a weak Hopf algebra
isomorphism $\varphi:H_2^1 \rar H$ with the following action;
\bea
\varphi(g^2) & = & I, \nn\\
\varphi(g)   & = & G, \nn\\ 
\varphi(x)   & = & X, \nn\\ 
\varphi(gx)  & = & GX.\nn
\eea

In a similar way to both of the examples above,
for a quantised enveloping algebra $U\equiv\Uqn$, its weak extension $U_w$ and some
other algebraic structure $U_0$, a decomposition of the form
$$
U_w = U_0 \oplus U = (1-J)U_w \oplus JU_w.
$$
exists due to there being a central idempotent $J$ whose existence
derives from the relaxation of the invertibility of group-like elements in
the algebra. In fact it is straightforward to prove the fact that
for any $d$, 
$$
\Uqn \simeq J.w\sl^d_q(n).
$$
This result is another way of stating Proposition 1 from the
paper \cite{LD}.

From another point of view, we could say that a weak
extension is nothing but the original Hopf algebra plus some
other algebra in which is contained all the information
regarding the weak structure. In this case, it would help to know
what conditions $U_0$ would have to satisfy in order that $U_w$
has a weak Hopf structure. This has not been the approach of this
paper as we saw in the last section. Since there are no other
weak Hopf algebra isomorphisms on the weak $w\sl^d_q(n)$, the
classification {\it as weak Hopf algebras} is complete. However, if we
were to consider all possible algebra isomorphisms, this direct
sum decomposition is important since it leads to discovering
several ``unexpected'' isomorphisms which do not arise from the
action of the automorphisms in the quantised enveloping algebra
case.

Since we only need the presence of a central idempotent to
achieve this direct sum decomposition, we can apply this idea to
the weak extensions of $\Uqn$ from the previous section, since
the element $J$ is always a central idempotent. However, it does
not affect our classification of the weak Hopf algebra structure from the
previous section. To demonstrate, 
we show that, rather unexpectedly, there is an
algebra isomorphism $\psi:w\sl^{10}_q(3)\rar w\sl^9_q(3)$. 

We first apply the direct sum decomposition to $U=w\sl^{10}_q(3)$
and $V=w\sl^9_q(3)$ such that
$$
U=(1-J)U\oplus JU \equiv U_0\oplus U_1
$$
and similarly
$$
V=(1-J')V\oplus J'V \equiv V_0\oplus V_1
$$
Explicitly we have $U_0$ generated by $<(1-J)E_1,(1-J)F_1,1-J>$
and $U_1$ generated by $<J,JE_1,JF_1,E_2,F_2,K_1,K_2,\kbar_1,\kbar_2>$.
For $V$, denoting its generators by a prime, we have $V_0$
generated by $<(1-J')E_1',(1-J')F_2',1-J'>$ and $V_1$ generated by 
$<J',J'E_1',E_2',F_1',J'F_2',K_1',K_2',\kbar_1',\kbar_2'>$. It is
straightforward to show that both $U_1$ and $V_1$ are isomorphic
as (weak) Hopf algebras to $U_q[sl(3)]$. It is also easy to verify that
$U_0$ and $V_0$ are both Abelian with the same number of
generators and are therefore isomorphic as algebras. Combining
these two facts leads to the isomorphism $\psi:U\rar V$, the
action of which is given by
\bea
\psi(1)   & = & 1, \nn\\
\psi(E_1) & = & (1-J')F_2'+F_1', \nn\\
\psi(E_2) & = & J'F_2', \nn\\
\psi(F_1) & = & E_1', \nn\\
\psi(F_2) & = & E_2', \nn\\
\psi(K_i) & = & \kbar_i', \nn\\
\psi(\kbar_i) & = & K_i'. \nn
\eea
The map $\psi$ is consistent with all the defining relations so
is therefore an algebra homomorphism and it can be shown to have
inverse
\bea
\psi^{-1}(E_1') & = & F_1, \nn\\
\psi^{-1}(E_2') & = & F_2, \nn\\
\psi^{-1}(F_1') & = & JE_1, \nn\\
\psi^{-1}(F_2') & = & (1-J)E_1+E_2, \nn\\
\psi^{-1}(K_i') & = & \kbar_i, \nn\\
\psi^{-1}(\kbar_i') & = & K_i. \nn
\eea
Therefore $\psi$ is an isomorphism.

If we set the action of the co-product $\d$ for $U$ and allow the
freedom to choose the co-product $\d'$ of $V$ consistently with
$\psi$, then we end up having to compare the action of  
$(\psi\tp\psi)\circ\d$ with $\d'\circ\psi$. Applying both of
these maps to the generators will clearly give a non-coassociative
$\d'$. Therefore the $\psi$ only can be considered as an algebra
isomorphism.
However, it is uncertain whether or not this co-product $\d'$
would define a quasi-bialgebra \cite{Dri2}. This certainly raises some
interesting questions relating to whether or not these
isomorphisms could correspond to some kind of Drinfeld twist. If
so, then perhaps our classification is only a much smaller classification 
of the structures as quasi-bialgebras. This idea may warrant further 
investigation.

%
%
%
%
\setcounter{equation}{0}
\section{Concluding remarks}

We have seen that it is possible to define weak extensions of
$\Uqn$ by only relaxing some of the relations in the original
algebra. As we saw in the work of Li and Duplij \cite{LD}, one
nice way of doing this is to relax invertibility of the
group-like elements to a more general regularity condition and also to 
impose one of two relations on the other generators. This allows
us to define many examples.

One observation is that it is also possible to extend the
definition of a quantised superalgebra \cite{K} to the weak case by 
using the same idea of relaxing the invertibility of the
generators $K$ and $\kbar$.
We demonstrate with the algebra $w\osp^d_q(2|1)$ which has generators $\{K,\kbar, V_+, V_-\}$.
We define the parity of these generators to be $p(K)=p(\kbar)=0$,
$p(V_\pm)=1$.

The following relations are satisfied;
\bea
  & & K \bK = \bK K, \qquad K \bK K = K, \qquad \bK K \bK = \bK, \nn \\
  & & K V_{\pm} = q^{\pm 1} V_{\pm} K, \qquad\qquad \bK V_{\pm} = q^{\mp 1}V_{\pm} \bK,
   \label{Wospdef1} \\
  & & \{ V_+, V_- \} = -\frac{1}{4} \frac{K-\bK}{\qqi}. \nn
\eea
Keeping in theme with the previous sections, if in addition 
the following relations are satisfied,
$$
K X \bK = q^{\pm 1} X
$$
where $X=V_+$ or $V_-$, then we call $X$ a type 2 generator.
Otherwise we call $X$ a type 1 generator.
This example is almost exactly like the case of $w\sl^d_q(2)$ in
that we have the same notion of generators of type 1 and 2. The
co-algebra structure is of the same form, and the only real
difference is that the weak antipode is a graded algebra anti-homomorphism,
so it satisfies $T(ab)=(-1)^{p(a)p(b)}T(b)T(a)$. 

All the weak Hopf algebras given in this paper 
have non-cocommutative co-products. 
This implies  existence of universal $R$-matrices that could give new solutions 
of quantum Yang-Baxter equations as mentioned in \cite{Li1,LD}.
One direction for future work is to investigate the form of such
$R$-matrices. We expect the expressions would not be that
different to those of the original Hopf algebra due to the direct
sum decomposition of section \ref{swee}. In fact, in section
\ref{swee} we gave one possible $R$-matrix for the finite
dimensional weak Hopf generalisation of Sweedler's well known
example using these facts.

Section \ref{swee} also demonstrates the fact that there
are many algebra isomorphisms between structures. In this
article we did not investigate all possible isomorphisms, but
instead gave a small existence proof that such isomorphisms do
indeed exist. It would be interesting to classify these
structures as algebras using this observation.

A question that arose during our investigation is one
related to automorphisms and twisting, especially in the usual quantised
enveloping algebra case. As we have already mentioned, the algebra 
automorphisms are well known for the quantised enveloping algebras, some
of which are also bialgebra automorphisms. We are currently
unaware of whether or not, corresponding to every algebra
automorphism $\psi: A\rar A$, there exists a twist element $F\in
A\tp A$ such that
$$
(\psi\tp\psi)\d(a) = F. \d(\psi(a)).F^{-1}, \ \ \forall a\in A.
$$

%
%
%
%

\end{document}